\documentclass[10pt,a4paper,reqno]{amsart}

\usepackage{amsmath,amssymb}
\usepackage{amsthm}
\usepackage{mathrsfs}
\usepackage{graphicx}
\usepackage{textcomp}
\usepackage{booktabs}

\newtheorem{thm}{Theorem}[section]

\newtheorem{lemma}[thm]{Lemma}

\theoremstyle{definition}

\newtheorem*{remark}{Remark}

\newtheorem*{prf}{Proof}

\numberwithin{equation}{section}

\def\amin{\alpha_{\text{min}}^{(k)}}
\def\aks{\alpha_{\ast}^{(k)}}
\def\akss{\alpha_{\ast\ast}^{(k)}}
\def\akp{\aks}
\def\akb{\aks}
\def\akt{\hat{t}^{(k)}}
\def\akc{\alpha_{\text{c}}^{(k)}}
\def\akh{\widehat{\alpha}^{(k)}}
\def\amax{\alpha_{\text{max}}^{(k)}}
\def\rk{\rho^{(k)}}
\def\rks{\rk_{\ast}}

\def\ek{\epsilon^{(k)}}
\def\lsk{\lambda_{\ast}^{(k)}}
\def\u#1{u_{1,#1}}
\def\PD#1#2{\partial_{#2}#1}

\def\q{\mathfrak{q}}
\def\qk{\q^{(k)}}
\def\Q{\mathfrak{Q}}
\def\Qk{\Q^{(k)}}
\def\QkN{\Q^{(k)}_N}
\def\QkD{\Q^{(k)}_D}
\def\h{\mathfrak{h}}
\def\hk{\h^{(k)}}
\def\hkone{\h_1^{(k)}}
\def\hktwo{\h_2^{(k)}}

\def\eig#1#2{\lambda_{#1,#2}}
\def\eigone#1{\eig{1}{#1}}
\def\eigtwo#1{\eig{2}{#1}}
\def\eigthree#1{\eig{3}{#1}}

\def\ed{\,\mathrm{d}}

\let\phi=\varphi
\let\epsilon=\varepsilon

\DeclareMathOperator{\spec}{Spec}

\title[Higher order anharmonic oscillators]%
{Spectral properties of higher order\\anharmonic oscillators}
\author{Bernard Helffer}
\author{Mikael Persson}
\address[Bernard Helffer]{D\'{e}partement de Math\'{e}matiques, 
B\^{a}timent 425, Univ Paris-Sud et CNRS, F-91405 Orsay
C\'{e}dex, France}
\email{Bernard.Helffer@math.u-psud.fr}
\address[Mikael Persson]{Aarhus University, Department of
Mathematical Sciences, 1530 Ny Munkegade, 8000 Aarhus C, Denmark}
\email{mickep@imf.au.dk}

\begin{document}
\begin{abstract}
We discuss spectral properties of the self-adjoint operator
\begin{equation*}
-\frac{d^2}{dt^2}+\Bigl(\frac{t^{k+1}}{k+1}-\alpha\Bigr)^2
\end{equation*}
in $L^2(\mathbb{R})$ for odd integers $k$. We prove that the minimum over 
$\alpha$ of the ground state energy of this operator is attained at a unique 
point which tends to zero as $k$ tends to infinity. Moreover, we show that
the minimum is non-degenerate. 
These questions arise naturally in the spectral analysis of Schr\"{o}dinger 
operators with magnetic field.
This extends or clarifies previous results
by Pan-Kwek~\cite{pakw}, Helffer-Morame~\cite{hemo1}, Aramaki~\cite{ara}, 
Helffer-Kordyukov~\cite{heko1,heko5,heko4} and Helffer~\cite{helf}.
\end{abstract}

\maketitle

\section{Introduction}

\subsection{Definition of $\Qk(\alpha)$ and main result}
For any $\alpha\in\mathbb{R}$ we denote by $\eigone{\Qk(\alpha)}$ the 
lowest eigenvalue  of the self-adjoint second order differential operator
\begin{equation*}
\Qk(\alpha)= -\frac{d^2}{dt^2}+\Bigl(\frac{t^{k+1}}{k+1}-\alpha\Bigr)^2.
\end{equation*}
We also denote by $\qk(\alpha)$ the quadratic form corresponding to 
$\Qk(\alpha)$,
\begin{equation*}
\qk(\alpha)[u]= \int_{\mathbb{R}}\bigl|u'(t)|^2
+\Bigl(\frac{t^{k+1}}{k+1}-\alpha\Bigr)^2|u(t)|^2 \ed t.
\end{equation*}
The main result of the present paper is the following theorem.

\begin{thm}\label{thm:main}
Assume that $k\geq 1$ is an odd integer. There exists a unique $\amin$ such that 
\begin{equation}
\inf_{\alpha\in\mathbb{R}} \eigone{\Qk(\alpha)}=\eigone{\Qk(\amin)}.
\end{equation}
Moreover, $\amin>0$ and the minimum is non-degenerate,
\begin{equation}
\partial^2_{\alpha\alpha}\bigl(\eigone{\Qk(\alpha)}\bigr)\bigm|_{\alpha=\amin}>0.
\end{equation}
\end{thm}

\begin{thm}\label{thm:even}
Assume that $k$ is even. Then $\alpha=0$ is a non-degenerate 
local minimum of $\eigone{\Qk(\alpha)}$.
\end{thm}

\begin{thm}\label{thm:largek}
If $k$ is odd, 
\begin{equation}\label{eq:aminzero}
\lim_{k\to +\infty} \amin = 0.
\end{equation}
In the even case, there exists $k_0$ such that for $k\geq k_0$ ($k$ even),
the ground state energy $\eigone{\Qk(\alpha)}$ has a unique minimum
which is attained at $\alpha=0$.
\end{thm}

\subsection{Historical context}
The operator $\Qk(\alpha)$ was first introduced in the context of magnetic 
Schr\"odinger operators in~\cite{mo}, and was further studied 
in~\cite{hemo1,pakw,heko1}.

The uniqueness of $\alpha_{min}^{(k)}$ was first observed numerically
 in~\cite{mo} for $k=1$. A proof for $k=1$ was given in~\cite{pakw}, 
which was completed in~\cite{helf}. The uniqueness for $k>1$ ($k$ odd) was 
announced in~\cite{ara} but the given proof seems incomplete. The non-degeneracy 
was obtained for $k=1$ in~\cite{helf} and conjectured in the general case
 in~\cite{heko5} and~\cite{heko4}. This conjecture was supported by numerical
computations performed by V.~Bonnaillie-No\"el, see Table~\ref{tab:bonn}. The 
results for large $k$ were announced in~\cite{heko5} and a proof was sketched 
in~\cite{heko6}.

\begin{table}[htbp]
\caption{Numerical values calculated by V.~Bonnaillie-No\"el with an accuracy 
of $10^{-2}$.}
\label{tab:bonn}
\begin{tabular}{crrrrrrrrrr}
\toprule
$k$ & 1& 2& 3& 4& 5& 6& 7& 8& 9& 10\\ 
\midrule
$\amin$ & 0.35& 0& 0.16& 0& 0.10& 0& 0.07& 0& 0.05& 0\\
$\eigone{\Qk(\amin)}$ & 0.57& 0.66& 0.68& 0.76& 0.81& 0.87& 0.92& 0.98& 1.02& 1.07\\
$\eigtwo{\Qk(\amin)}$ & 1.98& 2.50& 2.61& 2.98& 3.18& 3.47& 3.66& 3.90& 4.07& 4.27\\
$\eigthree{\Qk(\amin)}$ & 4.11& 5.24& 5.68& 6.52& 7.03& 7.69& 8.16& 8.70& 9.12& 9.57\\
\bottomrule
\end{tabular}
\end{table}

The outline of the paper is the following: In Section~\ref{sec:somefacts} we 
collect some facts about the operator $\Qk(\alpha)$, which we use 
in Section~\ref{sec:mainproof} to prove Theorem~\ref{thm:main}.
We prove Theorem~\ref{thm:even} in Section~\ref{sec:even}. We consider 
large values of $k$ in Section~\ref{sec:largek} and prove 
Theorem~\ref{thm:largek}.

\section{Auxiliary results}\label{sec:somefacts}

We recall some results about $\Qk(\alpha)$ obtained 
in~\cite{mo,heko5,heko4,heko6}.
\begin{lemma}\label{lem:alphabig}
It holds that $\eigone{\Qk(\alpha)}\to\infty$ as $|\alpha|\to\infty$.
\end{lemma}

\begin{prf}
We first note that if $k$ is odd and $\alpha<0$, then 
$\q^{(k)}(\alpha)[u]\geq \alpha^2\|u\|^2$, so 
$\eigone{\Qk(\alpha)}\geq \alpha^2$. On the other hand, for any integer $k>0$ 
one can use semi-classical analysis~\cite{si4,hesj} to show that
\begin{equation*}
\eigone{\Qk(\alpha)} \sim (k+1)^{2k/(k+1)}\alpha^{k/(k+1)},
\quad \alpha\to\infty.
\end{equation*}
For even $k$ it holds that $\eigone{\Qk(\alpha)}=\eigone{\Qk(-\alpha)}$.
\qed
\end{prf}

So, it is clear that the smooth function $\eigone{\Qk(\alpha)}$ is lower
semi-bounded, and
\begin{equation*}
\lsk := \inf_{\alpha\in\mathbb{R}} \eigone{\Qk(\alpha)}>0
\end{equation*}
and there exists (at least one) $\amin\in\mathbb{R}$ such that 
$\eigone{\Qk(\alpha)}$ is minimal,
\begin{equation*}
\eigone{\Qk(\amin)}=\lsk.
\end{equation*}
Let $\u{\alpha}\in L^2(\mathbb{R})$ be the $L^2$ normalized strictly positive
eigenfunction of the operator $\Qk(\alpha)$ corresponding to the eigenvalue
$\eigone{\Qk(\alpha)}$, 
\begin{equation}\label{eq:eq}
\Qk(\alpha)\u{\alpha} = \eigone{\Qk(\alpha)}\u{\alpha},\quad
\|\u{\alpha}\|=1.
\end{equation}
The function $\u{\alpha}$ can be chosen to depend smoothly on $\alpha$.

\begin{lemma}\label{lem:azero}
Assume that $k$ is odd. Then it holds that $\akc>0$ for all critical points 
$\akc$ of $\eigone{\Qk(\alpha)}$. In particular, $\amin>0$.
\end{lemma}
\begin{prf}
Differentiating~\eqref{eq:eq} with respect to $\alpha$ and taking the inner 
product with $\u{\alpha}$ we find
\begin{equation}\label{eq:lprim}
\partial_\alpha\eigone{\Qk(\alpha)}
 = -2\int_{-\infty}^\infty \Bigl(\frac{t^{k+1}}{k+1}-\alpha\Bigr)
\bigl(\u{\alpha}\bigr)^2 \ed t.
\end{equation}
So, when the derivative is zero, we get
\begin{equation}\label{eq:akceq}
\akc = \int_{-\infty}^\infty \frac{t^{k+1}}{k+1}
\bigl(\u{\akc}\bigr)^2 \ed t >0.
\end{equation}
\qed
\end{prf}

\begin{lemma}\label{lem:lthree}
Assume that $\akc$ is a critical point of $\eigone{\Qk(\alpha)}$. If either
\begin{itemize}
\item[(A)] $(k+2)\eigtwo{\Qk(\akc)}> (k+6)\eigone{\Qk(\akc)}$ or
\item[(B)] $k$ is odd or $\akc=0$, and 
$(k+2)\eigthree{\Qk(\akc)}> (k+6)\eigone{\Qk(\akc)}$,
\end{itemize}
then $\PD{}{\alpha\alpha}^2 \eigone{\Qk(\akc)} >0$. In particular this implies 
that $\eigone{\Qk(\alpha)}$ has a local minimum at $\akc$ which is 
non-degenerate.
\end{lemma}

\begin{prf}
We start by assuming that the condition in (A) is fulfilled.
The differentiation in the proof of Lemma~\ref{lem:azero} also provides us with 
a formula for $\partial_\alpha \u{\alpha}$,
\begin{equation*}
\partial_\alpha\u{\alpha} = -2\bigl(\Qk(\alpha)-\eigone{\Qk(\alpha)}\bigr)^{-1}
\Bigl[\Bigl(\frac{t^{k+1}}{k+1}-\alpha\Bigr)\u{\alpha}\Bigr],
\end{equation*}
where the inverse is the regularized resolvent. Differentiating~\eqref{eq:eq} 
twice, we find
\begin{equation*}
\partial^2_{\alpha\alpha}\eigone{\Qk(\alpha)} = 2-4\int_{-\infty}^\infty 
\Bigl(\frac{t^{k+1}}{k+1}-\alpha\Bigr)\u{\alpha}\partial_\alpha\u{\alpha}\ed t
\end{equation*}
By an application of the Cauchy-Schwarz inequality and the bound 
\begin{equation}\label{eq:resbound}
\|\bigl(\Qk(\alpha)-\eigone{\Qk(\alpha)}\bigr)^{-1}\|
\leq \frac{1}{\eigtwo{\Qk(\alpha)}-\eigone{\Qk(\alpha)}}
\end{equation}
we find that
\begin{equation}\label{eq:middle}
\partial^2_{\alpha\alpha}\eigone{\Qk(\alpha)} 
\geq 2-\frac{8}{\eigtwo{\Qk(\alpha)}-\eigone{\Qk(\alpha)}}
\Bigl\|\Bigl(\frac{t^{k+1}}{k+1}-\alpha\Bigr)\u{\alpha}\Bigr\|^2.
\end{equation}
To calculate the norm on the right-hand side, we note that the ground state
energy of the operator
\begin{equation*}
\Qk(\alpha,\rho)=-\frac{1}{\rho^2}\frac{d^2}{dt^2}
+\Bigl(\rho^{k+1}\frac{t^{k+1}}{k+1}-\alpha\Bigr)^2
\end{equation*}
is independent of $\rho$, i.e.,
\begin{equation*}
-\frac{1}{\rho^2}\frac{d^2}{dt^2}\u{\alpha,\rho}
+\Bigl(\rho^{k+1}\frac{t^{k+1}}{k+1}-\alpha\Bigr)^2\u{\alpha,\rho}
= \eigone{\Qk(\alpha)}\u{\alpha,\rho}.
\end{equation*}
Differentiating this identity with respect to $\rho$ and then letting $\rho=1$ 
and $\alpha=\akc$, and then taking the inner product with $\u{\akc}$, we 
get
\begin{equation*}
(k+1)\Bigl\|\Bigl(\frac{t^{k+1}}{k+1}-\akc\Bigr)\u{\akc}\Bigr\|^2 
= \Bigl\|\frac{d}{dt}\u{\akc}\Bigr\|^2,
\end{equation*}
and consequently
\begin{equation}\label{eq:akcpot}
\Bigl\|\Bigl(\frac{t^{k+1}}{k+1}-\akc\Bigr)\u{\akc}\Bigr\|^2
 =\frac{1}{k+2}\eigone{\Qk(\akc)}.
\end{equation}
Inserting this in~\eqref{eq:middle} we find that
\begin{align*}
\partial^2_{\alpha\alpha}\eigone{\Qk(\alpha)} 
&\geq 2-\frac{8\eigone{\akc}}
{(k+2)\Bigl(\eigtwo{\Qk(\akc)}-\eigone{\Qk(\akc)}\Bigr)}\\
&=2\frac{(k+2)\eigtwo{\Qk(\akc)}-(k+6)\eigone{\Qk(\akc)}}
{(k+2)\Bigl(\eigtwo{\Qk(\akc)}-\eigone{\Qk(\akc)}\Bigr)}.
\end{align*}
Hence, if for some $k$, we have
\begin{equation*}
(k+2)\eigtwo{\Qk(\akc)}>(k+6)\eigone{\Qk(\akc)}
\end{equation*}
we deduce that the minimum is non-degenerate. This finishes the proof under 
assumption (A). If, instead, (B) is satisfied, then we observe that 
$\partial_\alpha\u{\alpha}$ is an even function, and for even functions we 
have~\eqref{eq:resbound} with $\eigthree{\Qk(\alpha)}$ in place of 
$\eigtwo{\Qk(\alpha)}$. The rest follows the same lines as in the proof of (A).
\qed
\end{prf}

\begin{lemma}\label{lem:abound}
Assume that $k$ is odd and that $\akc$ is a critical point of 
$\eigone{\Qk(\alpha)}$. Then
\begin{equation}\label{eq:abound}
\bigl(\akc\bigr)^2 < \eigone{\Qk(\akc)}.
\end{equation}
\end{lemma}

\begin{prf}
Using the fact that $\u{\alpha}$ is even we get, using integration by parts,
\begin{equation*}
\int_0^\infty \frac{d}{dt}\Bigl[\Bigl(\frac{t^{k+1}}{k+1}-\alpha\Bigr)^2\Bigr]
\bigl(\u{\alpha}\bigr)^2 \ed t = \bigl(\eigone{\Qk(\alpha)}-\alpha^2\bigr)
\u{\alpha}(0)^2.
\end{equation*}
For a critical point $\alpha=\akc$, we get
\begin{equation*}
\int_0^\infty \Bigl(\frac{t^{k+1}}{k+1}-\alpha\Bigr) 
\bigl(\u{\akc}\bigr)^2 \ed t = 0.
\end{equation*}
Combining these two formulas, we obtain
\begin{multline*}
\Bigl(\eigone{\Qk(\akc)}-\bigl(\akc\bigr)^2\Bigr)
\u{\akc}(0)^2 \\
 =  2 \int_0^\infty \Bigl(t^k-\bigl(\akc(k+1)\bigr)^{k/(k+1)}\Bigl)
\Bigl(\frac{t^{k+1}}{k+1}-\akc\Bigr) 
\bigl(\u{\akc}\bigr)^2 \ed t > 0.
\end{multline*}
If $\u{\akc}(0)=0$, then $\u{\akc}\equiv 0$ since $\u{\akc}'(0)=0$, and
so~\eqref{eq:abound} holds.
\qed
\end{prf}

\section{Proof of Theorem~\ref{thm:main}}\label{sec:mainproof}
We will use the lemmas in the previous section to complete the proof. For that,
we need an upper bound on $\eigone{\Qk(\alpha)}$ and a lower bound on
$\eigthree{\Qk(\alpha)}$.
\subsection{Upper bound}
In this section we are looking for a good upper bound of $\eigone{\Qk(\alpha)}$.
\begin{lemma}\label{lem:upper}
For all $k\geq 1$ and $\alpha>0$ it holds that
\begin{equation}\label{eq:upper}
\eigone{\Qk(\alpha)}\leq \alpha^2 + 
\frac{\pi^2}{4}\frac{k+2}{k+1}
\Bigl(\frac{1}{4}(k+1)(2k+3)(2k+4)(2k+5)\Bigr)^{-1/(k+2)}
\end{equation}
In particular, if $k$ is odd, it holds that $\amin\leq\aks$ where
\begin{equation}\label{eq:aks}
\aks = \frac{\pi}{2}\Bigl(\frac{k+2}{k+1}\Bigr)^{1/2}
\Bigl(\frac{1}{4}(k+1)(2k+3)(2k+4)(2k+5)\Bigr)^{-1/(2k+4)}
\end{equation}
\end{lemma}

\begin{prf}
We will motivate our choice of trial function, inspired by~\cite{heko6}. 
For large $k$, the potential $\bigl(\frac{t^{k+1}}{k+1}-\alpha\bigr)^2$ will look
more and more as potential $p_{\alpha,\infty}$, 
\begin{equation*}
p_{\alpha,\infty}(t)=
\begin{cases}
\alpha^2 & |t|\leq 1\\
\infty & |t|>1.
\end{cases}
\end{equation*}
Among the potentials $p_{\alpha,\infty}$, $p_{0,\infty}$ is the one that will 
give the lowest energy, corresponding to the Dirichlet problem of 
$-\frac{d^2}{dt^2}$ on $L^2((-1,1))$, with eigenvalues 
\begin{equation*}
\Big\{\Bigl(\frac{\pi j}{2}\Bigr)^2\Bigr\}_{j\in\mathbb{N}\setminus\{0\}},
\end{equation*}
and with first eigenfunction $\cos(\pi t/2)$. Motivated by this, we introduce 
a parameter $\rho>0$ and use as a trial function
\begin{equation*}
u(t)=
\begin{cases}
\cos\bigl(\frac{\pi t}{2\rho}\bigr) & |t|\leq \rho,\\
0 & |t|>\rho.\\
\end{cases}
\end{equation*}
This function does not belong to the domain of $\Qk(\alpha)$, but to the form
domain of $\qk(\alpha)$, which is enough to use the min-max principle. 
A simple calculation shows that if $k$ is odd then
\begin{multline*}
\eigone{\Qk(\alpha)}\leq \frac{\qk(\alpha)[u]}{\|u\|^2} = \alpha^2 -
2\alpha \frac{\rho^{k+1}}{k+1}
\biggl(\frac{1}{k+2}+I\Bigl(\frac{k+1}{2}\Bigr)\biggr)\\
+\frac{\rho^{2k+2}}{(k+1)^2}\biggl(\frac{1}{2k+3}+I(k+1)\biggr)+
\frac{\pi^2}{4\rho^2},
\end{multline*}
where $I(m)=\int_0^1 s^{2m}\cos(\pi s)\ed s\leq0$. By integration by parts 
we see that
\begin{equation*}
-\frac{1}{2m+1}\leq I(m)\leq -\frac{1}{2m+1}+\frac{\pi^2}{(2m+1)(2m+2)(2m+3)}.
\end{equation*}
If $k$ is even the coefficient in front of 
$\alpha$ is zero. In any case we get
\begin{equation}\label{eq:minen}
\eigone{\Qk(\alpha)}\leq \alpha^2 + 
\frac{\pi^2\rho^{2k+2}}{(k+1)^2(2k+3)(2k+4)(2k+5)}+\frac{\pi^2}{4\rho^2}
\end{equation}
The right-hand side above is clearly minimal for $\alpha=0$. A differentiation
in $\rho$ also shows that it is minimal for 
\begin{equation*}
\rho=\rks:=\Bigl[\frac14(k+1)(2k+3)(2k+4)(2k+5)\Bigr]^{1/(2k+4)},
\end{equation*}
and if we put $\rks$ 
into~\eqref{eq:minen} and simplify we obtain~\eqref{eq:upper}.

The second statement is an immediate consequence of Lemma~\ref{lem:abound}.
\qed
\end{prf}

\begin{remark}
It holds that $\lim_{k\to+\infty} \rks=1$, which is coherent with the
fact that for the limiting case the first eigenfunction corresponds to $\rho=1$.
\end{remark}

\subsection{Lower bound on $\eigthree{\Qk(\alpha)}$}

\begin{lemma}\label{lem:lower}
Assume that $k\geq 3$ is odd and that $0\leq \alpha\leq\aks$, where $\aks$ is 
the constant from~\eqref{eq:aks}. Then
\begin{equation}\label{eq:lower}
\frac{k+2}{k+6}\eigthree{\Qk(\alpha)}\geq 
\eigone{\Qk(\alpha)}
\end{equation}
\end{lemma}

\begin{prf}
We introduce the operator $\QkN(\alpha)$ as the self-adjoint operator in
$L^2(\mathbb{R}^+)$ acting as
\begin{equation*}
\QkN(\alpha) =-\frac{d^2}{dt^2} + \Bigl(\frac{t^{k+1}}{k+1}-\alpha\Bigr)^2 
\end{equation*}
and with a Neumann condition at $t=0$. Since it holds that
$\eigtwo{\QkN(\alpha)}=\eigthree{\Qk(\alpha)}$ we will work on the half-line 
with $\QkN(\alpha)$ instead of $\Qk(\alpha)$, and show the inequality
\begin{equation*}
\frac{k+2}{k+6}\eigtwo{\QkN(\alpha)}\geq 
\eigone{\Qk(\alpha)}.
\end{equation*}

We introduce constants $0<\ek<1$ and $\akh>0$, to be determined 
in~\eqref{eq:ek} and~\eqref{eq:akh} below. We also 
set 
\begin{equation*}
\akt=\bigl((k+1)\akh\bigr)^{1/(k+1)}.
\end{equation*} 
We claim that if $0<\alpha<\ek\akh<\akh$, then
\begin{equation}\label{eq:pot}
\Bigl(\frac{t^{k+1}}{k+1}-\alpha\Bigr)^2 \geq 
p(t):=
\begin{cases}
\frac{2k\bigl(1-\ek\bigr)}{k+1}\bigl(\akt\bigr)^{2k}\bigl(t-\akt\bigr)^2 
& t>\akt\\
0 & 0<t\leq\akt.
\end{cases}
\end{equation}
This is clear for $0<t\leq\akt$. For $t>\akt$, we note that the the function 
$\hat{p}(t)=\bigl(\frac{t^{k+1}}{k+1}-\alpha\bigr)^2-p(t)$ is positive at
$t=\akt$, has a positive derivative at $t=\akt$,
\begin{equation*}
\hat{p}'(\akt)=2\biggl(\frac{\bigl(\akt\bigr)^{k+1}}{k+1}-\alpha\biggr)(\akt)^k
> 2\akh(1-\ek)(\akt)^k>0
\end{equation*}
and that $\hat{p}$ is convex for $t>\akt$,
\begin{equation*}
\begin{aligned}
\hat{p}''(t) &= 2t^{2k}+2\Bigl(\frac{t^{k+1}}{k+1}-\alpha\Bigr)kt^{k-1}
-\frac{4k\bigl(1-\ek\bigr)}{k+1}\bigl(\akt\bigr)^{2k}\\
&>2\bigl(1-\ek\bigr)(\akt)^{2k} +2k(\akh-\alpha)(\akt)^{k-1}
-\frac{4k\bigl(1-\ek\bigr)}{k+1}\bigl(\akt\bigr)^{2k}\\
&>\frac{2k\bigl(1-\ek\bigr)}{k+1}\bigl(\akt\bigr)^{2k} 
+ 2k\bigl(1-\ek\bigr)\akh (\akt)^{k-1} 
- \frac{4k\bigl(1-\ek\bigr)}{k+1}\bigl(\akt\bigr)^{2k}\\
&=0.
\end{aligned}
\end{equation*}
Let us denote by $\hk$ the self-adjoint operator in $L^2(\mathbb{R}^+)$,
acting as
\begin{equation*}
\hk = -\frac{d^2}{dt^2}+p(t),
\end{equation*}
and with a Neumann condition at $t=0$. Next, we decompose our Hilbert space
$L^2(\mathbb{R}^+)$ as 
$L^2(\mathbb{R}^+)=L^2((0,\akt))\oplus L^2((\akt,\infty))$
and introduce two new operators $\hkone$ and $\hktwo$.

The first one, $\hkone$, is the self-adjoint 
operator in $L^2((0,\akt))$ acting as
\begin{equation*}
\hkone=-\frac{d^2}{dt^2}, \quad 0<t<\akt
\end{equation*}
with Neumann boundary conditions at $t=0$ and $t=\akt$. This operator has 
eigenvalues 
\begin{equation*}
\spec\bigl(\hkone\bigr)=
\biggl\{\biggl(\frac{(j-1) \pi}{\akt}\biggr)^2\biggr\}_{j=1}^\infty.
\end{equation*}
The second operator, $\hktwo$, is the self-adjoint operator in 
$L^2((\akt,\infty))$, acting as
\begin{equation*}
\hktwo=-\frac{d^2}{dt^2}+\frac{2k(1-\ek)}{k+1}\bigl(\akt\bigr)^{2k}
\bigl(t-\akt\bigr)^2,\quad t>\akt 
\end{equation*}
with Neumann condition at $t=\akt$. After translation $s=t-\akt$ we get
\begin{equation*}
-\frac{d^2}{ds^2}+ \frac{2k(1-\ek)}{k+1}\bigl(\akt\bigr)^{2k}s^2,\quad s>0 
\end{equation*}
with Neumann condition at $s=0$. We use a scaling argument and compare with the 
harmonic oscillator on the half-line. The result is that the eigenvalues of 
$\hktwo$ are
\begin{equation*}
\spec\bigl(\hktwo\bigr)=
\biggl\{\biggl[\frac{2k(1-\ek)}{k+1}
\biggr]^{1/2}\bigl(\akt\bigr)^k(4j-3)\biggr\}_{j\in\mathbb{N}\setminus\{0\}}
\end{equation*}

We clearly have
\begin{equation*}
\eig{j}{\QkN(\alpha)}\geq \eig{j}{\hk}\geq \eig{j}{\hkone\oplus\hktwo},
\quad j\in\mathbb{N}\setminus\{0\},
\end{equation*}
and $\spec\bigl(\hkone\oplus\hktwo\bigr)
=\spec\bigl(\hkone\bigr)\cup \spec\bigl(\hktwo\bigr)$.

Next, we choose $\akh$ so that the second eigenvalue of $\hkone$ agrees with the
first one of $\hktwo$, i.e.,
\begin{equation*}
\Bigl(\frac{\pi}{\akt}\Bigr)^2 
= \biggl[\frac{2k(1-\ek)}{k+1}\biggr]^{1/2}\bigl(\akt\bigr)^{k}.
\end{equation*}
This gives
\begin{equation}\label{eq:akh}
\akt= \biggl[\frac{\pi^4(k+1)}{2k(1-\ek)}\biggr]^{\frac{1}{2(k+2)}},\quad 
\akh = \frac{1}{k+1}
\biggl[\frac{\pi^4(k+1)}{2k(1-\ek)}\biggr]^{\frac{k+1}{2(k+2)}},
\end{equation}
and the lower bound of $\eigtwo{\QkN(\alpha)}$ becomes 
\begin{equation*}
\eigtwo{\QkN(\alpha)}\geq \pi^2
\biggl[\frac{2k(1-\ek)}{\pi^4(k+1)}\biggr]^{\frac{1}{k+2}}.
\end{equation*}
Next we want to choose $\ek$ in such a way that both 
\begin{equation}\label{eq:inone}
\ek\akh\geq\akp, 
\end{equation}
and
\begin{equation}\label{eq:intwo}
\frac{k+2}{k+6}\pi^2
\biggl[\frac{2k(1-\ek)}{\pi^4(k+1)}\biggr]^{1/(k+2)}
\geq\alpha^2+\bigl(\akp\bigr)^2,\quad 0<\alpha\leq\akp
\end{equation}
are satisfied. It is clearly enough to prove the last 
inequality for $\alpha=\akp$. We let $\ek$ be given by
\begin{equation}\label{eq:ek}
\ek = 1-\frac{2}{k(k+1)}.
\end{equation}
With this choice, $\akt$ and $\akh$ reads
\begin{equation}\label{eq:akhnew}
\akt= \biggl[\frac{\pi^2(k+1)}{2}\biggr]^{\frac{1}{k+2}},\quad 
\akh = \frac{\pi^2}{2}\biggl[\frac{2}{\pi^2(k+1)}\biggr]^{\frac{1}{k+2}},
\end{equation}
and the lower bound of $\eigtwo{\QkN(\alpha)}$ becomes 
\begin{equation*}
\eigtwo{\QkN(\alpha)}\geq \pi^2
\biggl[\frac{2}{\pi^2(k+1)}\biggr]^{\frac{2}{k+2}}.
\end{equation*}
We start with \eqref{eq:inone}. We claim that $\ek\akh$ is monotonically
increasing for $k\geq 3$. Indeed, both factors are positive, and $\ek$ is 
obviously increasing.
We differentiate the expression for $\akh$ and use the fact that for $k\geq 3$
\begin{equation*}
\log\bigl(\pi^2(k+1)/2\bigr)>2,
\end{equation*}
to conclude that
\begin{equation*}
\frac{d}{dk}\akh = \akh
\Biggl[\frac{(k+1)\log\bigl(\pi^2(k+1)/2\bigr)-(k+2)}{(k+2)^2(k+1)}\Biggr] > 0.
\end{equation*}
Moreover, $\ek\akh$ is equal to 
$2^{-11/5}\times 3^{-1}\times 5\pi^{8/5}$ for $k=3$. We bound the constants
$\aks$ from above as
\begin{equation}\label{eq:aminbound}
\aks \leq \frac{\pi}{2}\sqrt{\frac{5}{4}}
\end{equation}
for $k\geq 3$.
Hence,~\eqref{eq:inone} is a consequence of
\begin{equation*}
2.26\approx 2^{-11/5}\times 3^{-1}\times 5\pi^{8/5} 
> \frac{\pi}{2}\sqrt{\frac{5}{4}}\approx 1.76.
\end{equation*}
For inequality~\eqref{eq:intwo}, we note that both sides are positive, so we 
will show that $A_1(k)\geq 1$ for all $k\geq 3$ with
\begin{equation}\label{eq:mylhs}
A_1(k):=\frac{\frac{k+2}{k+6}\pi^2
\Bigl[\frac{2}{\pi^2(k+1)}\Bigr]^{\frac{2}{k+2}}}
{2\bigl(\akp\bigr)^2} 
= \frac{2(k+1)}{k+6}
\Bigl[\frac{(2k+3)(2k+4)(2k+5)}{\pi^4(k+1)}\Bigr]^{\frac{1}{k+2}}.
\end{equation}
A plot of $A_1(k)$ is given in Figure~\ref{fig:Aone}. 
\begin{figure}[!h]
\begin{center}
\includegraphics{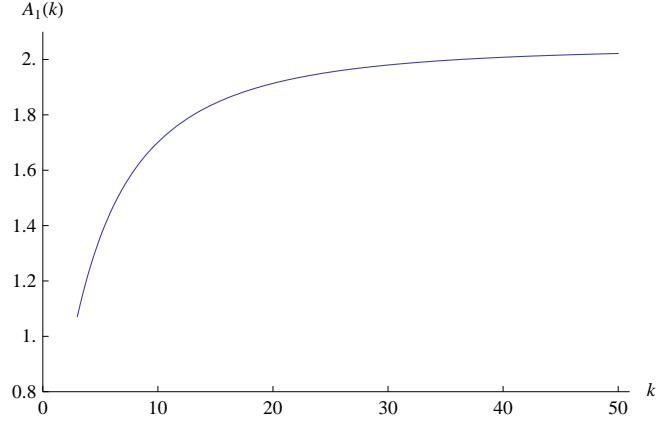}
\end{center}
\caption{A plot of $A_1(k)$ for $3\leq k\leq 50$}
\label{fig:Aone}
\end{figure}
Next, we use the estimate 
\begin{equation*}
(2k+3)(2k+4)(2k+5)>8(k+1)^3,
\end{equation*}
which implies that
\begin{equation*}
A_1(k)>\frac{2(k+1)}{k+6}
\Bigl(\frac{8(k+1)^2}{\pi^4}\Bigr)^{1/(k+2)}
\end{equation*}
The first factor is greater than $1$ if $k\geq 5$ and the second one is greater
than $1$ if $k\geq 3$. For $k=3$ get
\begin{equation*}
A_1(3)=2^{14/5}3^{-8/5}5^{1/5}11^{1/5}\pi^{-4/5}\approx 1.07.
\end{equation*}
This finishes the proof of~\eqref{eq:intwo} and completes the proof.
\qed
\end{prf}

\subsection{End of proof of Theorem~\ref{thm:main}}
By Lemmas~\ref{lem:azero} and~\ref{lem:abound} it follows that
\begin{equation*}
0<\amin<\akp.
\end{equation*}
However, by Lemmas~\ref{lem:lthree} and~\ref{lem:lower}, we find that all 
critical points in this interval must be non-degenerate minima. This
clearly implies the uniqueness, and finishes the proof.\qed

\section{The case of even $k$}\label{sec:even}

In this section we prove Theorem~\ref{thm:even}. 

\begin{prf}[of Theorem~\ref{thm:even}]
The lower bound of 
$\eigtwo{\Qk(\alpha)}$ from Lemma~\ref{lem:evenlower} is no good for small 
values of $\alpha$. Instead, we use the lower bound
\begin{equation*}
\Bigl(\frac{t^{k+1}}{k+1}-\alpha\Bigr)^2 
\geq \Bigl(\frac{|t|^{k+1}}{k+1}-\alpha\Bigr)^2,
\end{equation*}
and then we use that the second eigenvalue corresponding to the potential on
the right-hand side on $\mathbb{R}$ is equal to the first eigenvalue of
the operator
\begin{equation*}
\QkD(\alpha) = -\frac{d^2}{dt^2}+\Bigl(\frac{t^{k+1}}{k+1}-\alpha\Bigr)^2
\end{equation*}
in $L^2(\mathbb{R}^+)$ with a Dirichlet condition at $t=0$. We use the same type
of splitting as in Lemma~\ref{lem:lower},
\begin{equation*}
\QkD(\alpha)\geq \hk=\hkone\oplus\hktwo,
\end{equation*}
and write
\begin{equation*}
0\leq \alpha < \ek\akh < \akh,\quad \akt=(\akh(k+1))^{1/(k+1)},
\end{equation*}
where the constants $\ek$, $\akh$ and $\akt$ play the same roles as in the proof 
of Lemma~\ref{lem:lower} (but, as we will see, they are not the same!). 
This time the operator $\hkone$ is given by
\begin{equation*}
\hkone = -\frac{d^2}{dt^2}
\end{equation*}
in $L^2((0,\akt))$ with Dirichlet condition at $t=0$ and Neumann condition at
$t=\akt$. This operator has eigenvalues
\begin{equation*}
\spec\bigl(\hkone\bigr)=\biggl\{
\biggl(\frac{(2j-1)\pi}{2\akt}\biggr)^2
\biggr\}_{j\in\mathbb{N}\setminus\{0\}}.
\end{equation*}
The operator $\hktwo$ is the same as in the proof of Lemma~\ref{lem:lower}, 
with eigenvalues
\begin{equation*}
\spec\bigl(\hktwo\bigr)=
\biggl\{\biggl[\frac{2k(1-\ek)}{k+1}\bigl(\akt\bigr)^{2k}
\biggr]^{1/2}(4j-3)\biggr\}_{j\in\mathbb{N}\setminus\{0\}}
\end{equation*}

As in Lemma~\ref{lem:lower}, the best lower bound we can get on 
$\eigone{\QkD(\alpha)}$  is the one we get when the first eigenvalues of 
$\hkone$ and $\hktwo$ are equal. This determines $\akh$ as
\begin{equation}\label{eq:akheven}
\akt= \biggl[\frac{\pi^4(k+1)}{32k(1-\ek)}\biggr]^{\frac{1}{2(k+2)}}
,\quad 
\akh = \frac{1}{k+1}
\biggl[\frac{\pi^4(k+1)}{32k(1-\ek)}\biggr]^{\frac{k+1}{2(k+2)}},
\end{equation}
We let $\ek=\frac{1}{2k}$. Then the lower bound becomes
\begin{equation*}
\eigtwo{\QkN(\alpha)}\geq \frac{\pi^2}{4}
\biggl[\frac{32(k-1/2)}{\pi^4(k+1)}\biggr]^{\frac{1}{k+2}}.
\end{equation*}
To get the existence of an $\amax>0$ such that condition (A) in 
Lemma~\ref{lem:lthree} is fulfilled for $\alpha\in(-\amax,\amax)$ it is 
by Lemma~\ref{lem:upper} enough to show that $A_2(k)>1$ with
\begin{equation*}
A_2(k):=\frac{k+2}{k+6}\frac{\frac{\pi^2}{4}
\Bigl[\frac{32(k-1/2)}{\pi^4(k+1)}\Bigr]^{\frac{1}{k+2}}}
{\bigl(\akb\bigr)^2}
=\frac{k+1}{k+6}
\biggl(\frac{8}{\pi^4}(k-1/2)(2k+3)(2k+4)(2k+5)\biggr)^{1/(k+2)}.
\end{equation*}
See Figure~\ref{fig:Aktwo} for a plot of $A_2(k)$ for $2\leq k\leq 50$. We
note that $\lim_{k\to\infty}A_2(k)=1$. 
\begin{figure}[!h]
\begin{center}
\includegraphics{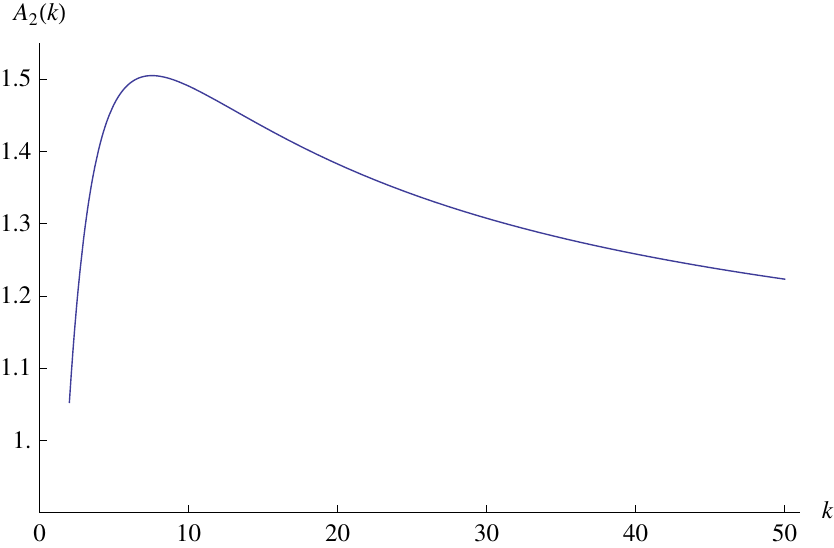}
\end{center}
\caption{A plot of $A_2(k)$ for $2\leq k\leq 50$}
\label{fig:Aktwo}
\end{figure}
By using the estimate 
\begin{equation*}
(k-1/2)(2k+3)(2k+4)(2k+5)>16(k+1)^3
\end{equation*}
(which is valid for all $k\geq 2$) we find that $A_2(k)>B(k)$ with
\begin{equation*}
B(k):=\frac{k+1}{k+6}\biggl(\frac{128}{\pi^4}(k+1)^3\biggr)^{1/(k+2)}.
\end{equation*}
The derivative of $B(k)$ is given by
\begin{equation*}
B'(k)=B(k)\frac{8k^2+44k+56-(k+1)(k+6)\log\bigl(\frac{128}{\pi^4}(k+1)^3\bigr)}
{(k+1)(k+2)^2(k+6)}.
\end{equation*}
For $k\geq 14$ it holds that $\log \bigl(\frac{128}{\pi^3} (k+1)^3\bigr)>8$ and 
so
\begin{equation*}
8k^2+44k+56-(k+1)(k+6)\log\Bigl(\frac{128}{\pi^4}(k+1)^3\Bigr)< 8-12k\leq -160,
\end{equation*}
which implies that $B'(k)<0$. Moreover, since $B(14)\approx 1.27$ and 
$\lim_{k\to\infty}B(k)=1$ it follows that $B(k)\geq1$, and thus $A_2(k)>1$, for 
all $k\geq 14$.

For even $2\leq k\leq 12$, we calculate $A_2(k)$ numerically,
\medskip
\begin{center}
\begin{tabular}{crrrrrr}
\toprule
$k$ & 2& 4& 6& 8& 10& 12\\ 
\midrule
$A_2(k)$ & 1.05& 1.41& 1.49& 1.50& 1.49& 1.47\\
\bottomrule
\end{tabular}
\end{center}
\medskip
which establishes $A_2(k)>1$ for all even $k\geq 2$.

The proof of the theorem is completed by an application of 
Lemma~\ref{lem:lthree}, noting that $\alpha=0$ is a critical point 
of $\eigone{\Qk(\alpha)}$ since $\eigone{\Qk(\alpha)}$ is even.
\qed
\end{prf}

\section{The case of large $k$}\label{sec:largek}

The goal of this section is to prove Theorem~\ref{thm:largek}. It will be done
using the ideas from~\cite{heko6}. 

For even $k\geq 2$ we introduce 
\begin{equation}\label{eq:mk}
m_k=\inf_{t\in\mathbb{R}}\frac{t^{k+1}-1}{t-1}.
\end{equation}
The constants $m_k$ decrease from $3/4$ for $k=2$ to $1/2$ as $k\to\infty$.

\begin{lemma}\label{lem:evenlower}
Let $k\geq 2$ be an even integer. With $m_k$ as in~\eqref{eq:mk} it holds that
\begin{equation}\label{eq:evenlower}
\eig{j}{\Qk(\alpha)}\geq m_k\alpha^{k/(k+1)}(k+1)^{-1/(k+1)}(2j-1),
\quad j\in\mathbb{N}\setminus\{0\}.
\end{equation}
\end{lemma}

\begin{prf}
We use a lower bound of the potential
\begin{multline*}
\Bigl(\frac{t^{k+1}}{k+1}-\alpha\Bigr)^2 
= \alpha^2
\biggl(\Bigl(\frac{t}{(\alpha(k+1))^{1/(k+1)}}\Bigr)^{k+1}-1\biggr)^2\\
\geq \Bigl(m_k \alpha^{k/(k+1)}(k+1)^{-1/(k+1)}\Bigr)^2
\Bigl(t-(\alpha(k+1))^{1/(k+1)}\Bigr)^2,
\end{multline*}
and then estimate with the eigenvalues of the harmonic oscillator on the whole 
line.
\qed
\end{prf}

\begin{lemma}\label{lem:akss}
Assume that $k\geq 2$ is an even integer and that $m_k$ is the constant
from Lemma~\ref{lem:evenlower}. Then $\amin\leq\akss$ where
\begin{equation}\label{eq:akss}
\akss\leq  
\biggl[\frac{(k+1)^{1/(k+1)}}{m_k}
\biggl(\frac{\pi^2}{4}\frac{k+2}{k+1}
\Bigl(\frac{1}{4}(k+1)(2k+3)(2k+4)(2k+5)\Bigr)^{-1/(k+2)}\biggr)
\biggr]^{(k+1)/k}.
\end{equation}
In particular, if $\eta>\frac{\pi^2}{2}$ then there exists $k_0$ such that, 
for $k\geq k_0$, $k$ even, $\eigone{\Qk(\alpha)}$ attains its minimum in 
$(-\eta,\eta)$.
\end{lemma}

\begin{prf}
Inequality~\eqref{eq:akss} follows by combining Lemma~\ref{lem:upper} 
(with $\alpha=0$) with Lemma~\ref{lem:evenlower}. The second statement is
immediate, by letting $k\to\infty$, and using the fact that $m_k\geq \frac12$ 
for all $k$.
\qed
\end{prf}

\begin{lemma}\label{lem:eiglargek}
Let $\alpha>0$. For any $j\in\mathbb{N}\setminus\{0\}$ it holds that
\begin{equation}\label{eq:eiglargek}
\lim_{k\to\infty}\eig{j}{\Qk(\alpha)} 
= \alpha^2 + \Bigl(\frac{j\pi}{2}\Bigr)^2
\end{equation}
with a uniform control with respect to $\alpha$ in any compact interval.
\end{lemma}

This result might be a consequence of $\Gamma$-convergence of the Pisa school,
except possibly for the uniform control of $\alpha$. See also~\cite{si5}, 
in particular Example~4.2. For the sake of completeness, we give a proof 
inspired by the methods in~\cite{fl}. 

\begin{prf}
We start with the upper bound, which we prove for $j\leq 2$ only. 
The general proof uses the same argument.

For $j=1$ the upper bound follows from Lemma~\ref{lem:upper}. For $j=2$, 
let us consider the functions
\begin{equation}
\phi_1(t)=
\begin{cases}
\cos\bigl(\frac{\pi t}{2}\bigr) & \text{if $|t|\leq 1$}\\
0 & \text{if $|t|>1$},
\end{cases}
\quad \text{and}\quad
\phi_2(t)=
\begin{cases}
\sin(\pi t) & \text{if $|t|\leq 1$}\\
0 & \text{if $|t|>1$}.
\end{cases}
\end{equation}
They are eigenfunctions of the two lowest eigenvalues of the limiting model 
$k\to\infty$, $-\frac{d^2}{dt^2}+\alpha^2$ in $L^2((-1,1))$ with Dirichlet 
boundary conditions.

Computing the energy of the function $\mu_1\phi_1+\mu_2\phi_2$, 
$|\mu_1|^2+|\mu_2|^2=1$, we find a sphere in a two-dimensional space on which 
the energy is less than $\mu(k)$, with
\begin{equation*}
\mu(k) = \alpha^2 + \pi^2 + C\frac{1+|\alpha|}{k+1}.
\end{equation*}
The upper bound in~\eqref{eq:eiglargek} for $j=2$ is a consequence of the 
min-max principle. We continue with the lower bound.

Let $\epsilon>0$ be given. Then, for bounded $\alpha>0$, we can choose $k$ so 
large that
\begin{equation*}
\Bigl(\frac{t^{k+1}}{k+1}-\alpha\Bigr)^2 \geq p(t):=
\begin{cases}
\bigl(\frac{(1+\epsilon)^{k+1}}{k+1}+\alpha\bigr)^2, 
& \hphantom{-1-\epsilon}\llap{$-\infty$}<t\leq -1-\epsilon,\\
\alpha^2(1-\epsilon), & -1-\epsilon < t \leq 1-\epsilon,\\
0, & \hphantom{-}1-\epsilon < t \leq 1+\epsilon,\\
\bigl(\frac{(1+\epsilon)^{k+1}}{k+1}-\alpha\bigr)^2, 
& \hphantom{-}1+\epsilon<t<\infty.
\end{cases}
\end{equation*}
We want to solve the eigenvalue equation
\begin{equation}\label{eq:plambda}
-\frac{d^2}{dt^2}u+p(t)u = \lambda u,
\end{equation}
by solving it for each interval and glue the solutions together as is done in
several examples in~\cite{fl}. We first note that the operator is positive, 
so we only have to consider $\lambda\geq 0$. Let us introduce the notation
\begin{gather*}
A = \Bigl(\frac{(1+\epsilon)^{k+1}}{k+1}+\alpha\Bigr)^2,\quad
B = \Bigl(\frac{(1+\epsilon)^{k+1}}{k+1}-\alpha\Bigr)^2,\quad
C = \sqrt{\strut\lambda-\alpha^2(1-\epsilon)},\\
t_0=-1-\epsilon,\quad t_1=1-\epsilon,\quad t_2=1+\epsilon.
\end{gather*}
We may choose $k$ so large that $A>\lambda$ and $B>\lambda$. 

If $\lambda>\alpha^2(1-\epsilon)$, the square integrable solution 
to~\eqref{eq:plambda} is given by
\begin{equation}
u(t)=
\begin{cases}
a_0\exp\bigl(\sqrt{A-\lambda}t\bigr) & -\infty<t\leq t_0,\\
b_0\cos(Ct)+b_1\sin(Ct) & \hphantom{-\infty}\llap{$t_0$} < t \leq t_1,\\
c_0\cos(\sqrt{\lambda}t)+c_1\sin(\sqrt{\lambda}t), 
& \hphantom{-\infty}\llap{$t_1$} < t \leq t_2,\\
d_0\exp\bigl(-\sqrt{B-\lambda}t\bigr), &\hphantom{-\infty}\llap{$t_2$}<t<\infty.
\end{cases}
\end{equation}
Here $a_0$, $b_0$, $b_1$, $c_0$, $c_1$ and $d_0$ are constants that are 
determined by gluing the solution together. 
The conditions that both $u$ and $u'$ 
should coincide at the points $t_0$, $t_1$ and $t_2$ read
\begin{equation*}
\begin{aligned}
a_0\exp\bigl(\sqrt{A-\lambda}t_0) &= b_0\cos(Ct_0)+b_1\sin(Ct_0),\\
a_0\sqrt{A-\lambda}\exp\bigl(\sqrt{A-\lambda}t_0)
&=-b_0C\sin(Ct_0)+b_1C\cos(Ct_0),\\
b_0\cos(Ct_1)+b_1\sin(Ct_1)
&=c_0\cos(\sqrt{\lambda}t_1)+c_1\sin(\sqrt{\lambda}t_1),\\
-b_0C\sin(Ct_1)+b_1C\cos(Ct_1)
&=-c_0\sqrt{\lambda}\sin(\sqrt{\lambda}t_1)
+c_1\sqrt{\lambda}\cos(\sqrt{\lambda}t_1),\\
c_0\cos(\sqrt{\lambda}t_2)+c_1\sin(\sqrt{\lambda}t_2)
&=d_0\exp\bigl(-\sqrt{B-\lambda}t_2\bigr),\\
-c_0\sqrt{\lambda}\sin(\sqrt{\lambda}t_2)
+c_1\sqrt{\lambda}\cos(\sqrt{\lambda}t_2)
&=-d_0\sqrt{B-\lambda}\exp\bigl(-\sqrt{B-\lambda}t_2\bigr).
\end{aligned}
\end{equation*}
This is a linear system of equations in $a_0$, $b_0$, $b_1$, $c_0$, $c_1$ and
$d_0$ which has nontrivial solutions if and only if 
\begin{equation}\label{eq:glue}
\frac{1}{C}
\frac{\sqrt{A-\lambda}\tan\bigl(C(t_1-t_0)\bigr)+C}
{\sqrt{A-\lambda}-C\tan\bigl(C(t_1-t_0)\bigr)}
=
-\frac{1}{\lambda}
\frac{\sqrt{B-\lambda}\tan\bigl(\sqrt{\lambda}(t_2-t_1)\bigr)+\sqrt{\lambda}}
{\sqrt{B-\lambda}-\sqrt{\lambda}\tan\bigl(\sqrt{\lambda}(t_2-t_1)\bigr)}
\end{equation}
This is the equation that determines the eigenvalues $\lambda$.
For large $k$, the terms $\sqrt{A-\lambda}$ and $\sqrt{B-\lambda}$ are 
dominating, and we can write~\eqref{eq:glue} as
\begin{equation}
\frac{1}{C}\tan\bigl(C(t_1-t_0)\bigr)
= -\frac{1}{\sqrt{\lambda}}\tan\bigl(\sqrt{\lambda}(t_2-t_1)\bigr)
+\mathcal{O}\bigl((k+1)(1+\epsilon)^{-(k+1)}\bigr)
\end{equation}
as $k\to\infty$, where the estimate is uniform for bounded $\alpha$ and 
$\lambda$. Inserting the values for $t_0$, $t_1$, $t_2$ and $C$, we find that
\begin{equation}\label{eq:coeffeq}
\frac{1}{\sqrt{\lambda-\alpha^2(1-\epsilon)}}
\tan\bigl(2\sqrt{\lambda-\alpha^2(1-\epsilon)}\bigr)
= -\frac{1}{\sqrt{\lambda}}\tan\bigl(2\epsilon\sqrt{\lambda}\bigr)
+\mathcal{O}\bigl((k+1)(1+\epsilon)^{-(k+1)}\bigr).
\end{equation}
If $0<\lambda<\alpha^2(1-\epsilon)$ then hyperbolic 
functions appear in the solution of~\eqref{eq:plambda}, and the same type of 
calculations that resulted in~\eqref{eq:coeffeq} this time yield
\begin{equation*}
\frac{1}{\sqrt{\alpha^2(1-\epsilon)-\lambda}}
\tanh\bigl(2\sqrt{\alpha^2(1-\epsilon)-\lambda}\bigr)
= -\frac{1}{\sqrt{\lambda}}\tan\bigl(2\epsilon\sqrt{\lambda}\bigr)
+\mathcal{O}\bigl((k+1)(1+\epsilon)^{-(k+1)}\bigr).
\end{equation*}
The function
\begin{equation*}
f_1(\lambda)= 
\begin{cases}
\frac{1}{\sqrt{\alpha^2(1-\epsilon)-\lambda}}
\tanh\bigl(2\sqrt{\alpha^2(1-\epsilon)-\lambda}\bigr),
& 0<\lambda<\alpha^2(1-\epsilon),\\
2, & \lambda=\alpha^2(1-\epsilon),\\
\frac{1}{\sqrt{\lambda-\alpha^2(1-\epsilon)}}
\tan\bigl(2\sqrt{\lambda-\alpha^2(1-\epsilon)}\bigr),
& \alpha^2(1-\epsilon)<\lambda<\infty,
\end{cases}
\end{equation*}
is positive for all 
$0\leq \lambda<\alpha^2(1-\epsilon)+\frac{\pi^2}{16}$, and
$\lim_{\lambda\nearrow \alpha^2(1-\epsilon)+\frac{\pi^2}{16}}f_1(\lambda) 
= +\infty$. 
For larger $\lambda$ it holds that $f_1(\lambda)$
is monotonically increasing from $-\infty$ to $+\infty$ in every interval
\begin{equation}\label{eq:interval}
\biggl(\biggl[\frac{(j-1/2)\pi}{2}\biggr]^2+\alpha^2(1-\epsilon),
\biggl[\frac{(j+1/2)\pi}{2}\biggr]^2+\alpha^2(1-\epsilon)\biggr).
\end{equation}
The function
\begin{equation*}
f_2(\lambda)= -\frac{1}{\sqrt{\lambda}}\tan\bigl(2\epsilon\sqrt{\lambda}\bigr)
\end{equation*}
is negative for all  $0\leq\lambda<\bigl(\frac{\pi}{4\epsilon}\bigr)^2$, and
$\lim_{\lambda\nearrow\bigl(\frac{\pi}{4\epsilon}\bigr)^2} f_2(\lambda)
= -\infty$.

We find that if $\epsilon$ satisfies
\begin{equation*}
\Bigl(\frac{\pi}{4\epsilon}\Bigr)^2
>\biggl[\frac{(j+1/2)\pi}{2}\biggr]^2+\alpha^2(1-\epsilon)
\end{equation*}
then there exists a $k_j(\epsilon)$ and $C_j$ such that if 
$k\geq k(\epsilon)$, $k$ even, 
it holds that the $j$th solution of~\eqref{eq:coeffeq} lies in the 
interval~\eqref{eq:interval} and we conclude that 
\begin{equation*}
\eig{j}{\Qk(\alpha)} 
\geq \biggl[\frac{(j-1/2)\pi}{2}\biggr]^2 +\alpha^2(1-\epsilon)-C_j\epsilon.
\end{equation*}
This is not the upper bound we wanted. However we can do better. There exists a 
constant $K_j>0$ (uniform in $\alpha$, $\epsilon$) such that
\begin{equation*}
0<\lambda<\biggl[\frac{(j+1/2)\pi}{2}\biggr]^2+\alpha^2(1-\epsilon)
\implies-K_j\epsilon < f_2(\lambda)<0.
\end{equation*}
This implies that the first $j$ solutions to~\eqref{eq:coeffeq}, up to an error
of order $\epsilon$ coincide with the first $j$ zeros of the function 
$f_1(\lambda)$, i.e., for all $\epsilon>0$ there exist 
$\hat{k}_j(\epsilon)$ and $\widehat{C}_j$ such that for 
$k\geq \hat{k}_j(\epsilon)$, $k$ even, it holds that
\begin{equation*}
\eig{j}{\Qk(\alpha)}
\geq \Bigl(\frac{j\pi}{2}\Bigr)^2 +\alpha^2(1-\epsilon)-\widehat{C}_j\epsilon.
\end{equation*}
This completes the proof of~\eqref{eq:eiglargek}.
\qed
\end{prf}

We are now ready to prove Theorem~\ref{thm:largek}.

\begin{prf}[of Theorem~\ref{thm:largek}]
First, we show~\eqref{eq:aminzero}, where we consider odd $k$ only. We recall 
the bound~\eqref{eq:aminbound} on $\amin$, 
$0<\amin<\sqrt{\frac54}\frac{\pi}{2}$, and the formula~\eqref{eq:akceq} which is
valid for $\amin$, i.e.,
\begin{equation*}
\amin = \int_{\mathbb{R}} \frac{t^{k+1}}{k+1}\bigl(\u{\amin}\bigr)^2 \ed t.
\end{equation*}
It is enough to show that
\begin{equation*}
\lim_{k\to\infty} 
\int_1^\infty \frac{t^{k+1}}{k+1}\bigl(\u{\amin}\bigr)^2\ed t=0.
\end{equation*}
We first show that, for any $\epsilon>0$ it holds that
\begin{equation}\label{eq:aminfirst}
\lim_{k\to\infty} 
\int_{1+\epsilon}^\infty \frac{t^{k+1}}{k+1}\bigl(\u{\amin}\bigr)^2 \ed t=0.
\end{equation}
For any $k\geq 3$ and $0<\alpha<\sqrt{\frac54}\frac{\pi}{2}$ we use 
Lemma~\ref{lem:upper} to find
\begin{multline}\label{eq:helpinf}
\int_{1+\epsilon}^\infty \frac{t^{2(k+1)}}{(k+1)^2}\bigl(\u{\alpha}\bigr)^2\ed t
\leq 2\int_{1+\epsilon}^\infty \Bigl(\frac{t^{k+1}}{k+1}-\alpha\Bigr)^2
\bigl(\u{\alpha}\bigr)^2\ed t\\
+2\alpha^2\leq 2\eigone{\Qk(\alpha)}+2\alpha^2 < \frac{15}{8}\pi^2.
\end{multline}
In particular, we get
\begin{equation*}
\int_{1+\epsilon}^\infty \frac{t^{k+1}}{k+1}\bigl(\u{\amin}\bigr)^2 \ed t 
\leq (1+\epsilon)^{-(k+1)}(k+1)\frac{15}{8}\pi^2,
\end{equation*}
which establishes~\eqref{eq:aminfirst}. We write the remaining integral as
\begin{multline*}
\int_{1}^{1+\epsilon} \frac{t^{k+1}}{k+1}\bigl(\u{\amin}\bigr)^2 \ed t = \\
\int_{1}^{1+\epsilon} \Bigl(\frac{t^{k+1}}{k+1}-\amin\Bigr)
\bigl(\u{\amin}\bigr)^2 \ed t
 + \amin \int_{1}^{1+\epsilon}\bigl(\u{\amin}\bigr)^2 \ed t,
\end{multline*}
and apply the Cauchy-Schwarz inequality and use~\eqref{eq:akcpot} to conclude 
that the first integral tends to zero as $k\to\infty$. For the second integral
we use the general inequality
\begin{equation*}
\int_a^b u(t)^2\ed t \leq 4 \int_{\frac{a+b}{2}}^b u(t)^2 \ed t 
+ 2(b-a)^2\int_a^b u'(t)^2 \ed t
\end{equation*}
with $a=1$ and $b=1+\epsilon$. We use~\eqref{eq:helpinf} to find that
\begin{multline*}
\int_{1+\frac{\epsilon}{2}}^{1+\epsilon} \bigl(\u{\amin}\bigr)^2\ed t 
\leq \Bigl(1+\frac{\epsilon}{2}\Bigr)^{-2(k+1)} (k+1)^2 
\int_{1+\frac{\epsilon}{2}}^{\infty} 
\frac{t^{2(k+1)}}{(k+1)^2}\bigl(\u{\amin}\bigr)^2\ed t\\
\leq \frac{15}{8}\pi^2 \Bigl(1+\frac{\epsilon}{2}\Bigr)^{-2(k+1)} (k+1)^2.
\end{multline*}
Moreover we use the inequality
\begin{equation*}
\int_{1}^{1+\epsilon} \bigl(\u{\amin}'\bigr)^2 \ed t \leq \eigone{\Qk(\amin)}
\leq \frac{5}{16}\pi^2,
\end{equation*}
to get, finally,
\begin{equation*}
\int_{1}^{1+\epsilon}\bigl(\u{\amin}\bigr)^2 \ed t \leq 
\frac{15}{2}\pi^2 \Bigl(1+\frac{\epsilon}{2}\Bigr)^{-2(k+1)} (k+1)^2
+\frac{5}{8}\pi^2\epsilon^2.
\end{equation*}
This achieves the proof of~\eqref{eq:aminzero}. 

We continue with the proof of the second statement. We know that $\alpha=0$ is a
non-degenerate local minima. By Lemma~\ref{lem:akss} it is enough to show that 
there exists a $k_0$ such that condition~(A) in Lemma~\ref{lem:lthree} holds for
all $k\geq k_0$, $k$ even, and all $0<\alpha\leq \eta$, where 
$\eta>\frac{\pi^2}{2}$ is arbitary. However, it is clear by 
Lemma~\ref{lem:eiglargek} that this can be done. 
\qed
\end{prf}

\section*{Acknowledgements}
The authors thank Yuri Kordyukov for many discussions and for allowing us to 
reproduce some proofs of~\cite{heko6}. We also thank S\o{}ren Fournais and 
Xingbin Pan for fruitful discussions. 
This work was started when the authors were at the Erwin Schr\"odinger 
Institute (ESI) in Vienna which is gratefully acknowledged.
MP is supported
by the Lundbeck foundation and by European Research Council under the
European Community's Seventh Framework Program (FP7/2007--2013)/ERC grant
agreement $202859$.

\bibliographystyle{abbrv}
\def\cprime{$'$}

\end{document}